\documentclass[a4paper,12pt,twoside]{article}
\usepackage[latin1]{inputenc}
\usepackage[T1]{fontenc}
\usepackage{amssymb}
\usepackage{dochier}
\usepackage[final]{lpretex}
\usepackage{title}
\usepackage{a4}
\usepackage{amscd}
\usepackage{caption}
\input xy
\xyoption{all}

\EndOfDump

\def\DHpartformat#1{(#1) }
\def\DHrefpart#1{(\DHRefpart{#1})}
\LBsetstyleof{partno}{arabic}

\let\define\def

  \def\P {{\mathbb P}} 
 \def\R {{\mathbb R}}

\def\Z {{\mathbb Z}} 

\define \n {\mathbb N}
\define \z {\mathbb Z}
\define \q {\mathbb Q}
\define \PP {\mathbb P}

\def\sA {{\Cal A}}  
 \def\sE {{\Cal E}} \def\sF {{\Cal F}}
  \def\sI {{\Cal I}}
  \def\sL {{\Cal L}}
 \def\sN {{\Cal N}} \def\sO {{\Cal O}}
 
\def\sS {{\Cal S}}  \def\sU {{\Cal U}}
  \def\sX {{\Cal X}}

\define \cN {\Cal N}
\define \cf {\Cal F}
\define \cg {\Cal G}
\define \cE {\Cal E}
\define \ce {\Cal E}
\define \cc {\Cal C}
\define \cV {\Cal V}
\define \cA {\Cal A}
\define \cK {\Cal K}
\define \cO {\Cal O}
\define \cF {\Cal F}
\define \cn {\Cal N}
\define \cI {\Cal I}
\define \sP {\Cal P}
\define \sW {\Cal W}
\def\tA {\widetilde{\Cal A}}

\def\a {\alpha} \def\b {\beta} \def\g {\gamma}

\define \x {\xi}
\define \y {\eta}
\define \G {\Gamma}
\define \r {\rho}
\define \w {\omega}

\def \tZ {\widetilde Z}
\def\tX {\widetilde X}

\def \tD {\widetilde D}

\def \trho {\widetilde {\rho}}

\def \tp {\widetilde{\mathbb P}}
\define \tH {\widetilde H}
\define \tG {\widetilde{\Gamma}}
\define \tW {\widetilde W}
\define \tF {\widetilde F}
\define \tm {\widetilde m}
\define \St {\widetilde S}
\define \Xt {\widetilde X}
\define \tS {\widetilde S}
\define \tpsi {\widetilde \psi}
\define \txi {\widetilde \xi}
\define \tL {\widetilde L}
\define \tE {\widetilde E}
\define \tl {\widetilde l}
\define \tA {\widetilde A}
\define \tom {\widetilde\omega}
\define \tT {\widetilde T}
\define \tB {\widetilde B}
\define \tf {\widetilde f}
\define \tsA {\widetilde{\sA}}

\define \tM {\widetilde M}
\define \tphi {\widetilde{\phi}}
\define \trho {\widetilde{\rho}}
\define \tR {\widetilde R}
\define \tp {\widetilde p}
\define \tq {\widetilde q}
\define \tc {\widetilde c}
\define \tsF {\widetilde {\sF}}
\define \tsN {\widetilde {\sN}}
\define \tsU {\widetilde {\sU}}
\define \tsX {\widetilde {\sX}}
\define \th {\widetilde h}

\def\pd {\partial}

\def \Dx1 {\frac{\pd}{{\pd} x_1}}
\def \Dy1 {\frac{\pd}{{\pd} y_1}}
\def \Dz1 {\frac{\pd}{{\pd} z_1}}
\def \Dx2 {\frac{\pd}{{\pd} x_2}}
\def \Dy2 {\frac{\pd}{{\pd} y_2}}
\def \Dz2 {\frac{\pd}{{\pd} z_2}}

\def\q {\quad} 

\def\mapdiagr#1{\Big\searrow\rlap{$\raise 5pt\vbox{{\hbox{$\mkern -15mu\scriptstyle#1$}}}$}}   

\def\mapdiagl#1{\llap{$\raise 5pt\vbox{{\hbox{$\scriptstyle#1\mkern
-15mu$}}}$}\Big\swarrow}              

\def\Mapdiagr#1{\nearrow\rlap{$\lower 5pt\vbox{{\hbox{$\mkern
-15mu\scriptstyle#1$}}}$}} 

\def\Mapdiagl#1{\llap{$\lower 5pt\vbox{{\hbox{$\scriptstyle#1\mkern
-15mu$}}}$}\searrow} 

\def\Mapswr#1{\swarrow\rlap{$\lower 5pt\vbox{{\hbox{$\mkern
-15mu\scriptstyle#1$}}}$}}              

\def\Mapnwl#1{\nwarrow\rlap{$\lower 5pt\vbox{{\hbox{$\mkern
-15mu\scriptstyle#1$}}}$}}

\def \inj {\hookrightarrow}

\define \Rhook {\hookrightarrow}

\def \half {\raise1pt\hbox{$\scriptstyle
        \frac{1}{2}\displaystyle$}}

\def \x{{\sl X}\llap{$\mkern -2mu {\scriptstyle -}$}}
\def \Pic {\operatorname{Pic}}
\def \Sing {\operatorname{Sing}}
\define \Kod {\operatorname{Kod}}
\define \dimension {\operatorname{dim}}
\define \codim {\operatorname{codim}}
\define \contr {\operatorname{contr}}
\define \rk {\operatorname{rank}}
\define \im {\operatorname{im}}
\define \Mor {\operatorname{Mor}}
\define \Cl {\operatorname{Cl}}
\define \Hilb {\operatorname{Hilb}}
\define \degree {\operatorname{deg}}
\define \mult {\operatorname{mult}}
\define \Aut {\operatorname{Aut}}
\define \NS {\operatorname{NS}}
\define\MW {\operatorname{MW}}
\define \Gal {\operatorname{Gal}}
\define \ch {\operatorname{char}}
\define \Jac {\operatorname{Jac}}
\define \Km {\operatorname{Km}}
\define \Sec {\operatorname{Sec}}
\define \Stab {\operatorname{Stab}}
\define \Br {\operatorname{Br}}
\define \inv {\operatorname{inv}}
\define \tr {\operatorname{tr}}
\define \Frob {\operatorname{Frob}}
\define \Symn {\operatorname{Sym}^n}
\define \Ev {\sE^\vee}
\define \ordp {\operatorname{ord}_p}
\define \Supp {\operatorname{Supp}}
\define \Ann {\operatorname{Ann}}
\define \disc {\operatorname{disc}}
\define \Lie {\operatorname{Lie}}
\define \embdim {\operatorname{embdim}}

\def \hZ{\widehat Z}

\def\Tors{\operatorname{Tors}}

\def\hod#1#2#3#4{\ensuremath{\if#30 H^{#2}({#1},{\cal O}_{#1}) \else 
 H^{#2}(#1,\Omega^{#3}\if\relax{#4}\relax_{#1}\else _{#1/#4}\fi)\fi}}

\begin{document}
\title{On exceptional Enriques surfaces}
\author{T. Ekedahl${}^\dagger$}
\author{N. I. Shepherd-Barron}
\address
{Math. Dept.\\
King's College London\\
Strand\\
WC2R 2LS\\
UK}
\email{N.I.Shepherd-Barron@kcl.ac.uk}
\maketitle
\begin{abstract}
  We give a complete description of all classical (``$\Z/2$'')
  Enriques surfaces with non-zero global vector fields.
  In particular we show that such surfaces exist.
  The result that we obtain also applies to supersingular
  (``$\a_2$'') surfaces that fulfil a rather special condition.
  In the course of the classification we study some properties
  of genus 1 fibrations special to characteristic $2$ as well as make
  a close study of the genus 1 fibrations on the surfaces in question.
\end{abstract}

\begin{section}{Introduction}
  Whether a classical Enriques surface has a non-zero global
  vector field is of interest in, for example, the deformation
  theory of Enriques surfaces. In \cite{SB96} it was claimed
  that this cannot occur. However, there is an error in the
  proof and in fact the truth is the opposite: such surfaces
  do exist. It turns out that a certain condition, which, in the case
  of a classical Enriques surface, is equivalent to the existence of a
  vector field, is also of interest in the non-classical case;
  we call surfaces that satisfy this condition \emph{exceptional}.
  (See \ref{0.3} for the precise definition.) Our main results
  are summarized in Theorems \ref{A} to \ref{C}, where the conductrix is a certain
  effective divisor supported on the image of the singular
  locus of the canonical double cover. In particular,
  we shall see that the exceptional surfaces
  are just those that possess a special genus $1$
  fibration with a double fibre of type $\tE_{6,7,8}$.
  Such surfaces have been classified by Salomonsson \cite{Sa03}.
  We shall say that an
  exceptional surface is \emph{of type $T$} if the
  support of the conductrix is a $T$-configuration.
  We also say that a genus $1$ fibration on an Enriques
  surface is \emph{special} if it has a $2$-section
  isomorphic to $\P^1$ (it always has a $2$-section of arithmetic
  genus at most $1$). Such a $2$-section is called special.
  (Classically surfaces with such a pencil are referred to as
  special, and in \cite{CD89} the terminology
  \emph{degenerate $U$-pair} is used.)

  If $\Delta$ is a simply laced affine Dynkin diagram
  then $F_\Delta$ will denote the corresponding Kodaira--N{\'e}ron
  fibre.

  If $A$ is a divisor whose support is a diagram
  of type $T_{p,q,r}$ then we write
  $$A=(a;b_1, ..., b_{p-1};c_1,..., c_{q-1};d_1,..., d_{r-1})$$
  in an obvious way. For example, $F_{\tE_r}$ is
  $(3;2,1;2,1;2,1), (4;2;3,2,1;3,2,1)$ and
  $(6;3;4,2;5,4,3,2,1)$ according as $r=6,7$ or $8$.

  \begin{theorem}\label{A} Suppose that
    $X$ is an Enriques surface in characteristic $2$.

    \part[i]\label{A1} $X$ is exceptional if and only if its conductrix
    $A$ is of type $T_{p,q,r}$, $(p,q,r)$
    is one of $(2,3,7),(2,4,5),(3,3,3)$  and $A$ is, accordingly,
    $$(5;2;3,2;4,4,3,3,2,1),\ (3;1;2,2,1;2,2,1,1),\ (2;1,1;1,1;1,1).$$

    \part[ii]\label{A3} $X$ is exceptional if it has a
    quasi--elliptic
    fibration with a simple $\tE_{7,8}$ fibre. $X$ is
    of type $T_{3,3,3}$ or $T_{2,4,5}$ accordingly.

    \part[iii]\label{A2} $X$ is exceptional if and only if it
    has a special genus $1$ fibration with a double fibre
    of type $\tE_{6,7,8}$. It is then of
    type $T_{3,3,3},\ T_{2,4,5}$ or $T_{2,3,7}$, respectively.
    \noproof
  \end{theorem}

  \ref{A1} is Proposition \ref{3.11} and
  Theorem \ref{3.14}, \ref{A2} is  Proposition \ref{Aiii} and
  \ref{A3} is Theorem \ref{descript}.
        
  The definition of an exceptional Enriques surface is given by
  a simple condition on the conductrix but we also give the
  following elaboration of that condition.

  \begin{theorem}\label{B} An exceptional Enriques surface $X$
    is either a $\Z/2$-surface or an $\a_2$-surface. A
    $\Z/2$-surface $X$ is exceptional precisely when it has global
    vector fields, and then $\dim H^0(X,T_X)=1$.
    An $\a_2$-surface is exceptional exactly when the
    cup product on $H^1(X,\sO_X)\times H^0(X,\Omega^1_X)$
    is zero.
    \noproof
  \end{theorem}

        This is proved in Propositions \ref{unipotent} and \ref{2.8}.

        \begin{remark} The presence of vector fields on a
          $\Z/2$-surface
          clearly makes its deformation theory ``pathological''.
          We shall show elsewhere that an Enriques surface is
          exceptional exactly when a versal deformation of
          it as a unipotent Enriques surface is singular.
        \end{remark}

        We go on to discuss the classification of exceptional
        surfaces and show that all three types exist, both
        for $\Z/2$-surfaces and for $\a_2$-surfaces. We also
        describe all genus $1$ fibrations on them. This description
        is complicated for surfaces $X$ of type $T_{3,3,3}$, where we need
        to distinguish between surfaces of different
        \emph{MW-rank} $\MW(X)$. By definition,
        $$\MW(X)=8-\sum_s n(s) -1,$$
        where $s$ runs over the fibres of the unique elliptic
        pencil on $X$ and $n(s)$ is the number of irreducible
        components of $s$. ($\MW(X)$
        equals the
        Mordell--Weil rank of the Jacobian surface
        associated to $X$.)

        \begin{theorem}\label{C}

          \part[i]\label{C1} An exceptional Enriques surface of type $T_{2,3,7}$
          has a unique genus $1$ fibration. This fibration is
          quasi--elliptic.

          \part[ii]\label{C2} An exceptional Enriques surface of type $T_{2,4,5}$
          has $2$ or $3$ genus $1$ fibrations. These fibrations are
          all quasi--elliptic.

          \part[iii]\label{C3} An exceptional Enriques surface $X$ of type $T_{3,3,3}$
          has a unique elliptic $1$ fibration. There exist
          quasi--elliptic fibrations on it; these fibrations
          are arranged in triples and the set of triples
          is a torsor under a discrete group $G$. The group $G$
          is trivial if $\MW(X)=0$, while $G=\Z$ if $\MW(X)=1$
          and $G$ is the Weyl group of type $\tA_2$
          if $\MW(X)=2$, the maximum possible value.
          Each quasi--elliptic fibration appears in $2$
          of these triples if $\MW(X)\ne 0$.
          \noproof
        \end{theorem}

        \ref{C1} is Theorem \ref{3.15}, \ref{C2} is Theorem \ref{4.4}
        and \ref{C3} is a combination of Theorems \ref{MW2},
        \ref{MW1} and \ref{MW0}.
        
        We also give, in Theorem \ref{4.4},
        a description of the $(-2)$ curves on
        an exceptional surface of type $(2,4,5)$ or $(2,3,7)$.
        For surfaces of type $(3,3,3)$ they are described
        by Theorems \ref{MW2}, \ref{MW1} and \ref{MW0}.

        The base field of all varieties appearing here will be,
        unless explicitly noted otherwise, algebraically
        closed and of characteristic $2$.

        We shall name the types of Enriques surfaces after
        the $\tau$ in their $Pic^\tau$. In other terminology
        $\mu_2$-surfaces are \emph{singular},
        $\Z/2$-surfaces are \emph{classical}
        and $\a_2$-surfaces are \emph{supersingular}.
        When
        $\tau=\Z/2$ or $\a_2$ we refer to the
        surface as \emph{unipotent}. After Proposition
        \ref{unipotent} all surfaces will, unless stated otherwise,
        be unipotent.

        We shall use the extended Dynkin diagram notation
        for the normal crossing singular fibres of a relatively
        minimal genus $1$
        fibration.
        The $E$-series of (extended) Dynkin diagrams
        are also graphs of type $T_{*,*,*}$
        (cf. \cite{CD89}, p. 105)
        and we shall pass freely between the two kinds of notation.
      \end{section}
      \begin{section}{Preliminaries}

        \begin{lemma}\label{0.1} Suppose that $S$
          is an affine Noetherian scheme and that
          $\pi:X\to S$ is a proper morphism of relative
          dimension $\le n$. Assume that $H^n(X,\sO_X)\ne 0$
          and that $H^n(Z,\sO_Z)=0$ for all closed subschemes
          $Z$ of $X$ such that $Z\ne X$. Then

          \part[i] $H^0(X,\sO_X)$ is a field and

          \part[ii] $\sO_X$ has no non-zero subsheaves whose support
          is of relative dimension less than $n$.
          \begin{proof} By assumption, $H^n(X,-)$ is right exact
            on quasi--coherent sheaves. Suppose that
            $0\ne\lambda \in R:= H^0(X,\sO_X)$
            and let $X_\lambda$ denote the closed subscheme
            of $X$ defined by $\lambda$. So there is an exact
            sequence
            $$\sO_X\stackrel{\lambda}{\to}\sO_X\to\sO_{X_\lambda}\to
            0$$
            whose cohomology gives an exact sequence
            $$H^n(X,\sO_X)\stackrel{\lambda}{\to}
            H^n(X,\sO_X)\to H^n(X_\lambda,\sO_{X_\lambda})
            \to 0.$$
            By our assumptions, $H^n(X_\lambda,\sO_{X_\lambda})=0$,
            so that multiplication by $\lambda$ is surjective
            on $H^n(X,\sO_X)$. Now $H^n(X,\sO_X)$ is a non-zero
            finitely
            generated $\Gamma(S,\sO_S)$-module, and
            so a finitely generated $R$-module. Therefore
            it has a non-zero quotient $M$ killed by some maximal ideal
            $\mathfrak m$ of $R$. If $0\ne\lambda\in\mathfrak m$
            then $M=\lambda M=0$, contradiction.
            So $\mathfrak m=0$ and $R$ is a field, as claimed.

            \DHrefpart{ii}: If $Z$ is the closed subscheme of
            $X$ defined by $\sI$ then $H^n(Z,\sO_Z)=H^n(X,\sO_X)\ne
            0$,
            while $H^n(Z,\sO_Z)=0$ by assumption.
          \end{proof}
        \end{lemma}

        
        We apply this lemma to a particular divisor on an
        Enriques surface.

        \begin{lemma}\label{0.2} Suppose that $D$ is an effective
          divisor on an Enriques surface $X$, that $h^0(X,\sO(D))=1$
          and that $h^1(\sO_D)\ne 0$. Then $D$ contains a half-fibre
          of a genus $1$ fibration.
          \begin{proof} Note that, by Riemann--Roch and the assumption
            that $h^0(X,\sO(D))=1$, $D$ contains no effective subdivisor $E$
            with $E^2>0$.

            By Lemma \ref{0.1} and the Noetherian property,
            there exists $0\ne E\subset D$ that is minimal for the
            condition that $h^1(E,\sO_E)\ne 0$. By Lemma \ref{0.1} again,
            $h^0(E,\sO_E)=1$. So $\chi(E,\sO_E)\le 0$ and then,
            by Riemann--Roch, $E^2\ge 0$. So $E^2=0$
            and now the result follows from \cite{CD89}, Th. 3.2.1.
          \end{proof}
        \end{lemma}
        
        Suppose that $X$ is a Gorenstein scheme and that
        $$\xymatrix{
          {\tZ}\ar[r]^\pi\ar[dr]_\b&{Z}\ar[d]^\a\\
          &{X}
        }
        $$
        is commutative. Suppose also that $\a,\b$ are finite and flat
        of degree $2$ and that $\pi$ is finite and birational.
        Suppose that $\sI_C\subset \sO_Z$ is the conductor of $\pi$;
        then $\sI_C$ is also an ideal in $\sO_{\tZ}$.

        \begin{lemma}\label{cond def} $\sI_C$ is an invertible
        $\sO_{\tZ}$-module and there is an effective
        Cartier divisor $A$ on $X$ such that $\sI_C=\b^*\sO_X(-A)$.
        \begin{proof} Note first that the $\sO_X$-modules
          $\sL=\sO_Z/\sO_X$ and $\sL'=\sO_{\tZ}/\sO_X$ are invertible
          and there is a commutative diagram
          $$\xymatrix{
            {0}\ar[r]&{\sO_X}\ar[r]\ar[d]&{\sO_Z}\ar[r]\ar[d]&{\sL}\ar[r]\ar[d]&{0}\\
            {0}\ar[r]&{\sO_X}\ar[r]&{\sO_{\tZ}}\ar[r]&{\sL'}\ar[r]&{0.}
          }
          $$
          So $\sL'=\sL\otimes\sO_X(A)$ for some effective
          Cartier divisor $A$ on $X$. By adjunction,
          $\omega_{Z/X}\cong \a^*\sL^{-1}$ and $\omega_{\tZ/X}\cong
          \b^*\sL'^{-1}$. Since also
          $\omega_{\tZ/X}=\sI_C\pi^*\omega_{Z/X}$ the result follows.
        \end{proof}
      \end{lemma}
                
        We shall call $A$ the \emph{conductrix} of the $X$-morphism
        $\pi$ and $B:=2A$ the \emph{biconductrix}.

        \begin{lemma}\label{foliation} Assume that $\a$ is inseparable, so that
          there is a factorization
          $$\tZ\stackrel{\b}{\to} X\stackrel{\g}{\to}\tZ^{(1)}$$
          of the Frobenius $\tZ\to\tZ^{(1)}$ and $\g:X\to \tZ^{(1)}$
          is the quotient of $X$ by a $2$-closed foliation $\sF$
          of rank $1$. Then $\sF\cong \omega_X(B)$,
          so that $c_1(\sF)$ is numerically equivalent to $B$.
          \begin{proof} There is an exact sequence
            $$0\to\sO_X\to \b_*\sO_{\tZ}\to \omega_X(A)\to 0.$$
            Recall that
            the map
            $$\b_*\sO_{\tZ}\to\Omega^1_X: f\mapsto df^2$$
            induces an injective map
            $\sO_X(B)\cong Frob_X^*(\omega_X(A))\to\Omega^1_X$
            which is saturated, since $\tZ$ is normal.
            That is, there is a short exact sequence
            \begin{eqnarray}\label{0.6}
              0\to \sO_X(B)\to\Omega^1_X\to\sI_W\omega_X(-B)\to 0.
              \end{eqnarray}
            The dual of this is
            $$0\to\sF\to T_X\to \sI_W(-B)\to 0.$$
          \end{proof}
        \end{lemma}

        This foliation $\sF$ is the \emph{natural foliation}
        on $X$.

        When $X$ is an Enriques surface, $Z$ its canonical double cover
        and $\tZ\to Z$ is the normalization we shall also speak of
        the conductrix and biconductrix of $X$. 

        If $f:X\to S$ is a genus $1$ fibration in characteristic $2$
        the we get a map $f':X'\to S$ which is the
        pullback of $f$ by the Frobenius on $S$,
        and we also have the normalization $\rho:\tX\to X'$.
        The conductrix of $f$ is then, by definition,
        the conductrix of $\rho$. This leads to slight ambiguity
        because there are two conductrices, one of the surface
        and one of the fibration. However, this should cause no
        confusion.

        \begin{definition}\label{0.3}
          An Enriques surface whose biconductrix is $B$
          is \emph{exceptional}
          if $H^1(B,\sO_B)\ne 0$.
        \end{definition}

        \begin{proposition}\label{unipotent} An Enriques surface
          which is not unipotent 
          has empty conductrix and so is
          not exceptional.
          \begin{proof} The canonical double cover of a
            surface that is not unipotent is {\'e}tale.
          \end{proof}
        \end{proposition}

        So from now on we shall only consider Enriques surfaces $X$
        that are unipotent. There is a diagram
        $$\xymatrix{
          {\hZ}\ar[r]^{\delta} &
          {\tZ}\ar[r]^{\pi}\ar[dr]_{\b}&
          {Z}\ar[d]^{\a}\\
          &&{X}
        }$$
        where $\a:Z\to X$ is the canonical double cover
        and is inseparable, $\tZ\to Z$
        is the normalization and $\hZ\to\tZ$ the minimal resolution.
        The conductrix is $A$ and the biconductrix is $B$.
        We shall assume that $A\ne 0$.

        \begin{proposition}\label{0.5}\label{2.8}
          
          \part[i] $B$ is the divisorial part of the zero locus of any
          global $1$-form.

          \part[ii] $h^0(B,\sO_B)=1$. In particular, $A$ cannot
          contain a fibre or half-fibre in any genus $1$ fibration on
          $X$ and $\Supp A$ is a normally crossing configuration $\G$
          of $(-2)$-curves.

          \part[iii]   $A$ is $1$-connected,
          $D^2<0$ for all effective $D\le A$, $A^2=-2$
          and $\G$ is a tree.

          \part[iv] $\tZ$ is rational and has either $4$ singularities of type $A_1$ or
          $1$ of type $D_4$.

          \part[v] $X$ is exceptional if and only if $B$ contains a
          half-fibre
          of some genus $1$ fibration.

          \part[vi] If $X$ is a $\Z/2$-surface then $X$ is exceptional
          if and only if it has a vector field. In any case
          $h^0(X,T_X)\le 1$.

          \part[vii] If $X$ is an $\a_2$-surface then it is
          exceptional
          if and only if the cup product
          $$H^1(X,\sO_X)\otimes H^0(X,\Omega^1_X)\to H^1(X,\Omega^1_X)$$
          is identically zero.
          \begin{proof} \DHrefpart{i} This follows at once
            from Lemma \ref{foliation}.
            
            \DHrefpart{ii} Computing Chern classes of the sheaves
            appearing in the exact sequence \ref{0.6}
            shows that $\deg W-B^2=12$.
            Since $h^{10}(X)=1$ we get $h^0(X,\sO_X(B))=1$,
            which is the first part of \DHrefpart{ii}.
            The rest follows at once.

            \DHrefpart{iii} By \cite{CD89}, Prop. 3.1.2 and Th. 3.2.1, we have
            $h^0(X,\sO_X(2D))\ge 2$ if $D$ is effective and $D^2\ge 0$,
            so that $D^2<0$ for all effective $D\le A$.

            Castelnuovo's criterion shows that
            $\hZ$ is rational, since $A>0$, and therefore $\tZ$
            has du Val singularities. Since also
            $$\chi(\sO_{\tZ}) = \chi(\sO_X)+\chi(\sO_X(A))=A^2/2 + 2$$
            we see that $A^2= -2$.

            Now suppose that $A=C+D$ with $C,D>0$. Then
            $$-2=A^2=C^2+D^2+2C.D\le -4 +2C.D,$$
            so that $C.D\ge 1$. Now it is enough to observe
            that $A$ cannot contain a cycle $D$ with $D^2\ge 0$
            since $h^0(X,\sO_X(B))=1$.

            \DHrefpart{iv} $B^2=-8$ and $\deg W=4$. If $P\in\Sing\tZ$
            then we can write
            $$\sO_{\tZ,P}\hat{}=k[[x,y,z]]/(z^2-f(x,y))$$
            where $f\in\mathfrak m_{(x,y)}^2$ and
            $4\ge\deg_PW=\dim_k k[[x,y]]/(f'_x,f'_y)$.
            Calculation shows that then
            $(\tZ,P)$ is of type either $A_1$ or $D_4$.
            Since $\sum_P\deg_PW=4$ the proof
            of \DHrefpart{iv} is complete.

            
            \DHrefpart{v} If $0<D'\le D$ are divisors on $X$
            such that $h^1(D,\sO_D)=0$, then $h^1(D',\sO_{D'})=0$.
            This gives one direction, and the other follows from
            Lemma \ref{0.2}. 

            \DHrefpart{vi} Notice that, by duality,
            $H^1(B,\sO_B)\ne 0$ if and only if
            $H^0(B,\omega_B)\ne 0$. Furthermore there
            is a short exact sequence
            \begin{eqnarray}\label{0.7}
              0\to\omega_X\to\omega_X(B)\to\omega_B\to 0.
            \end{eqnarray}
            In the $\Z/2$ case we exploit the dual of the sequence
            \ref{0.6}. From this it follows that, if $H^0(X,T_X)\ne
            0$,
            then $H^0(X,\omega_X(B))\ne 0$. However,
            $H^0(X,\omega_X)=0$ and we conclude by taking the
            cohomology
            of \ref{0.7}.
            
            \DHrefpart{vii} Assume that $X$ is an $\a_2$-surface.
            Note that
            $$h^1(X,\sO_X)= h^0(X,\Omega^1_X)=1.$$
            Suppose that $\b\in H^1(X,\sO_X)$
            and $\eta\in H^0(X,\Omega^1_X)$.
            By \ref{0.6} and the fact that $W\ne 0$, it follows
            that $\eta$ is the image of some
            $\eta'\in H^0(X,\omega_X(B))$.
            Then $\eta\beta$ is the image of $\eta'\beta$.
            Since $H^0(X,\sO_X(B))=1$ we can suppose that
            $\eta'$ comes from $1\in H^0(X,\sO_X)$
            under the inclusion of \ref{0.7}.
            Thus $\eta'\b$ is the image of $\b$
            under the map $H^1(X,\sO_X)\to H^1(X,\sO_X(B))$.
            Since $\omega_X$ is trivial the inclusion
            $\sO_X\inj\sO_X(B)$ is isomorphic
            to the inclusion $\omega_X\inj\omega_X(B)$.
            It follows from \ref{0.7} and the fact that
            $h^0(\omega_X)=h^0(\omega_X(B))=1$
            that the map
            $H^1(X,\omega_X)\to H^1(X,\omega_X(B))$
            is zero exactly when $h^0(\omega_B)\ne 0$,
            which, as we have noted,
            is equivalent to $X$ being exceptional.
          \end{proof}
        \end{proposition}

        \begin{lemma} If $X$ is exceptional then $\G$
          is not negative definite.
          \begin{proof} If $\G$ is negative definite
            then there is a contraction $h:X\to Y$ of $\G$
            to a du Val singularity. Since $R^1h_*\sO_X=0$, it follows
            that $H^1(B,\sO_B)=0$ for all divisors $B$ supported
            on $\Delta$.
          \end{proof}
        \end{lemma}

        \begin{lemma}\label{fibr} Suppose that $g:X\to\P^1$ is
          a genus $1$ fibration.
          
          \part[i] $\a:Z\to X$ factors through the pullback $X_F$
          of $g$ by the Frobenius on $\P^1$. The map
          $Z\to X_F$ is an isomorphism outside the double
          fibres of $g$.
          \part[ii] The restriction of $\a$ to a half-fibre
          of $g$ 
          is non-trivial.
          \begin{proof} We claim that the restriction of $\a$
            to a simple fibre $\Phi$ of $g$ is trivial.

            For this, suppose first that $X$ is an $\a_2$-surface.
            Since $H^0(\Phi,\sO_\Phi)=k$, the base field,
            it is enough to show that the map
            $H^1(X,\sO_X)\to H^1(\Phi,\sO_\Phi)$
            is zero. This follows from the cohomology
            of the short exact sequence
            $$0\to \sO_X(-\Phi)\to\sO_X\to\sO_\Phi\to 0$$
            and the fact that $h^1(X,\sO_X(-\Phi))=1$.

            In the $\Z/2$ case we need only observe that
            $K_X$ is the difference of the two half-fibres
            so its restriction to $\Phi$ is obviously trivial.

            Consider now the sheaf $\sA=g_*\a_*\sO_Z$ on $\P^1$.
            Its restriction to $\Phi$ is a trivial vector bundle,
            so that $\sA$ is of rank $2$.

            Now assume that $X$ is an $\a_2$-surface.
            Then there is a
            short exact sequence
            $$0\to \sO_X\to \a_*\sO_Z \to \sO_X\to 0$$
            whose cohomology gives an exact sequence
            $$0\to\sO_{\P^1}\to\sA\to\sO_{\P^1}\to R^1g_*\sO_X.$$
            Since $\rk\sA=2$ the image of the boundary map is torsion.
            Now $\Tors R^1g_*\sO_X$ has length $1$ so
            $\sA/\sO_{\P^1}\cong\sO_{\P^1}$ or $\sO_{\P^1}(-1)$.
            Since $H^0(\P^1,\sO_{\P^1})=k$ is follows that
            $\sA/\sO_{\P^1}\cong\sO_{\P^1}(-1)$.
            Therefore $\Sp\sA$ is obtained by taking the
            square root of a non-zero quadratic form.

            The same statement holds for $\Z/2$-surfaces,
            and is easier to prove as $R^1g_*\sO_X$ is torsion-free
            and $g_*\omega_X=\sO_{\P^1}$.

            Up to isomorphism there are only two such covers of
            $\P^1$,
            the trivial one and the Frobenius map
            $F:\P^1\to\P^1$. Now $\sA$ is
            reduced
            so the cover is non-trivial, so is the Frobenius.
            So $Z\to\P^1$ factors through $F:\P^1\to \P^1$.
            Finally, the map $g^*\sA \to \a_*\sO_Z$
            is an isomorphism away from the double fibres.

            As for the second part, the case of a $\Z/2$-surface
            is well known. So suppose that $X$ is an $\a_2$-surface
            and that $\Psi$ is the half-fibre.
            Consider the exact sequence
            $$0\to\sO_X(-\Psi)\to\sO_X\to\sO_\Psi\to 0;$$
            this time $h^1(\sO_X(-\Psi))=0$ and so the map
            $H^1(X,\sO_X)\to H^1(\Psi,\sO_\Psi)$ is injective,
            which is what is needed.
          \end{proof}
        \end{lemma}
        \begin{corollary}\label{fibr2}\label{2.12} Suppose that
          $g:X\to\P^1$ is a genus $1$ fibration.

          \part[i] The natural foliation $\sF$ is the kernel of the derivative
          $g_*:T_X\to g^*T_{\P^1}$.

          \part[ii] $g:X\to\P^1$ factors through the quotient map
          $X\to X/\sF=\tZ^{(1)}$.

          \part[iii] If $g$ is quasi--elliptic and $R=R_g$ is its curve
          of cusps then $R\subset A$ and
          $R$ is a $(-2)$-curve on $X$.
          \begin{proof} \DHrefpart{i} and \DHrefpart{ii} are restatements of part (1)
            of Lemma \ref{fibr} and \DHrefpart{iii} is an immediate
            consequence.
          \end{proof}
        \end{corollary}
        For the rest of this paper $R_g$ will denote the curve of cusps
        in a quasi--elliptic fibration $g:X\to\P^1$.
      \end{section}
      \begin{section}{Genus $1$ fibrations and the conductrix}
      \label{admiss}
        Suppose that $X$ is a unipotent Enriques surface whose
        conductrix $A$ is non-zero. We know that $\G:=\Supp A$ is
        a tree of $(-2)$-curves.

        Fix a genus $1$ fibration $f:X\to S=\P^1$. For $s\in S$
        we write $f^{-1}(s)=d_sX_s$, where $d_s=d=1$ or $2$ is the
        multiplicity of the fibre and $X_s$ is a Kodaira--N{\'e}ron
        divisor. We can write
        $A=R_1+\sum A_s$ where each $A_s$ is supported on
        $X_s$, $R_1$ is horizontal (that is, finite over $S$).
        Moreover, $R_1=0$ if $f$ is elliptic and $R_1=R_f$
        if $f$ is quasi--elliptic.
        
        \begin{lemma}\label{easy ell} Assume that $\G$ is
          not negative definite.

          \part[i] $f$ is elliptic if and only if
          $\G$ is of type $\tD$ or $\tE$.

          \part[ii] If $f$ is elliptic
          then $\G$ is the support of a unique 
          $X_s$.

          \part[iii] $X$ has at most one elliptic fibration.
          \begin{proof} Since $f$ is generically smooth, $R_1=0$.
            \DHrefpart{i} then follows from the $1$-connectedness
            and the non-negativity of $\G$.

            \DHrefpart{ii} and \DHrefpart{iii} are immediate.
          \end{proof}
        \end{lemma}

        \begin{lemma}\label{transverse}\label{3.2}
          If $E$ is a $(-2)$-curve on $X$
          then $A.E\le 1$. If $E_1,E_2$ are $(-2)$-curves
          that meet transversely in at least one point then
          $A.(E_1+E_2)\le 0$.
          \begin{proof} 
            Consider the natural foliation $\sF\cong\omega_X(2A)$.
            If $E$ is generically tangent to $\sF$ then
            $c_1(\sF).E\le 2$; conversely, if $c_1(\sF).E=2$
            then $\sF$ is everywhere tangent to $E$.
            If $E$ is not generically tangent to $\sF$ then
            $c_1(\sF).E\le -2$.
            If $E_1$ and $E_2$ are transverse at a point
            then it is impossible to have both $c_1(\sF).E_1=2$
            and $c_1(\sF).E_2=2$, and the lemma is proved.
          \end{proof}
        \end{lemma}
        
        \begin{proposition} $\G$ is a chain $A_n$ or a tree $T_{p,q,r}$.
          \begin{proof}
            We need to show that $\G$ has no vertex of valency
            at least $4$ and does not have two vertices of valency
            $3$; that is, 
            that $\G$ contains no configuration of type
            $\tD_{n\ge 4}$.

            If $\G\supset\tD_4$ then,
            as an immediate consequence of Lemma \ref{transverse},
            $A\ge F_{\tD_4}$. But this contradicts Proposition
            \ref{2.8}.

            Similarly, if $\G\supset \tD_{n\ge 5}$ then, by \ref{2.8},
            there is a curve $C$ in the spine of $\tD_n$ of
            multiplicity $1$ in $A$. This leads quickly to
            a contradiction to Lemma \ref{transverse}.
          \end{proof}
        \end{proposition}

        \begin{lemma} $n\le 11$ and $p+q+r\le 12$, respectively.
          \begin{proof} $\rho(X)=10$.
          \end{proof}
        \end{lemma}

        \begin{lemma} If $\G$ is not negative definite then
          there is an affine Dynkin diagram $\Delta\subset\G$
          and for every such $\Delta$, $\G-\Delta$ is connected.
          \begin{proof} The existence is clear.

            For any such $\Delta$, there is a genus $1$ fibration
            $f:X\to\P^1$ such that $\Delta$ is the support of a fibre.
            Since $R_f$ is irreducible and each connected component
            of $\G-\Delta$ contains a vertex corresponding to
            a component of $R_f$, the lemma is
            proved.
          \end{proof}
        \end{lemma}

        {\bf{From now on we assume that $\G$ is not negative definite.}}
        So $\G$ is one of $T_{3,3,r\ge 3}$, $T_{2,4,s\ge 4}$ and
        $T_{2,3,t\ge 6}$.

        \begin{lemma} $r\le 4$, $s\le 5$ and
          $t\le 7$.
          \begin{proof} Suppose $\G=T_{3,3,r\ge 5}=
            T_{3,3,3}\cup (R,E_1)$ or
            $T_{3,3,3}\cup(R,E_1,E_2)$. Then there
            is a quasi--elliptic fibration $f:X\to\P^1$
            such that $R=R_f$ and there is a singular fibre
            $X_s$ such that $X_s\cap\Supp A\supset E_1$.
            Since $\rho(X)=10$ the fibre $X_s$ has at most $3$
            components, so $X_s=E_1+E_2=Ta$ or
            $E_1+E_2+E_3=Tr$ or $I_3$.
            Then either $A.E_2\ge 2$ or $A.(E_2+E_3)\ge 2$,
            both of which are impossible.

            If $\G=T_{2,4,6}$ then
            $\G=T_{2,4,4}\cup (R,E_1)$ and we get a quasi--elliptic
            fibration $f:X\to\P^1$ with a singular fibre $X_s$
            that contains $E_1$. Then $X_s=Ta$ and we get
            a contradiction as before.
          \end{proof}
        \end{lemma}

        Suppose that $A'=(a(r);a(r-1),...,a(1);b,b';c)$
        is that part of $A$ supported on a subdiagram
        $T_{2,3,r}$ of $\G$. So
        $$a(i-1)-2a(i)+a(i+1)\le 1\ \textrm{and}\
        a(i-1)-a(i)-a(i+1)+a(i+2)\le 0$$
        for all $i$. Define $d(i)=a(i+1)-a(i)$,
        so that $d(i)\le d(i-2)$ and
        $d(i)\le d(i-1)+1.$

        Take the least $j$ such that $d(j)<0$, if such exists.
        Then $d(0),\ldots, d(j-1)\ge 0$
        and $d(j),\ldots, d(r-1)<0$.

        \begin{lemma}\label{lem 1} $a(i)$
          increases (not necessarily strictly) as $i$
          goes from $1$ to $j+1$ and then decreases
          (not necessarily strictly) as $i$ goes from $j+1$ to $r$.
          \begin{proof} This is immediate, from the
            properties of $d$ just observed.
          \end{proof}
        \end{lemma}

        \begin{lemma}\label{lem 2} $a(r-1)<a(r)$.
          \begin{proof} There is a diagram $\Delta\subset\G$
            of type $D_5=T_{2,2,3}$
            such that the part $A''$ of $A$ that is supported
            on $\Delta$ is $A''=(a(r);a(r-1);b,b';c)$ where the
            entries are all $>0$. Then
            $$a(r-1)-a(r)-c+b\le 0\ \textrm{and}$$
            $$a(r-1)+b'+c-b-a(r)\le 0.$$
            Adding these gives
            $$1\le b'\le 2(a(r)-a(r-1)).$$
          \end{proof}
        \end{lemma}

        Therefore $a(i)$ increases as $i$ goes from $1$ to $r$.

        \begin{lemma}\label{lem 3} If $a(r)\ge r$ then
          $a(i)\ge i$ for all $i\in[1,r]$.
          \begin{proof} We know that $d(i)\ge 0$ for all $i$.
            Suppose that $d(m)=0$ for some $m$;
            take $m$ to be minimal. Since $d(i)\le d(i-2),$
            it follows that
            $$0=d(m)=d(m+2)=\cdots = d(m+2s)$$
            for all $s$. Since also $d(k+1)\le d(k)+1$,
            it follows that
            $$d(m+1),\ldots , d(m+2s+1)\le 1.$$
            So $a(r-2k)\ge a(r)-k)$ and
            $a(r-2k-1)\ge a(r)-k-1,$
            so that $a(i)\ge i$ for all $i\ge m$.

            In the range $i<m$ we have $d(i)\ge 1$, so $a(i)\ge i$
            for all $i<m$.
          \end{proof}
        \end{lemma}

        \begin{corollary}\label{cor 4}
          Suppose that $\G$ is affine. Then $a(r)<r$.
          \begin{proof} Suppose $a(r)\ge r$.
            Then $A\ge F_{\G}$, by Lemma \ref{lem 3},
            which is impossible.
          \end{proof}
        \end{corollary}

        Corollary \ref{cor 4} and Lemmas \ref{lem 1} and \ref{lem 2}
        make it easy to enumerate the cycles $A$
        such that $h^0(X,2A)=1$, $\Supp A$ is
        an affine diagram $\G$ and A satisfies
        the conclusions of Lemma \ref{transverse}
        (but maybe $A^2\ne -2$).
        \begin{enumerate}
        \item $\G=T_{3,3,3}:\ A=(2;1,1;1,1;1,1),\ A^2=-2.$

        \item $\G=T_{2,4,4}:\ A=(3;2,2,1;2,2,1;1),\ A^2=-2.$

          \item $\G=T_{2,3,6}:$
      \begin{enumerate}
      \item $A=(5;4,4,3,3,2;3,2;1),\ A^2=-4.$
      \item $A=(5;4,4,3,3,1;3,2;2),\ A^2=-4.$
      \item $A=(5;4,4,3,2,1;3,2;2),\ A^2=-2.$
      \end{enumerate}
      \end{enumerate}
      If $\G=T_{p,q,r}$ is hyperbolic then $\G=\G_{aff}\cup\{R_g\}$
      where $R_g$ is the curve of cusps in a quasi--elliptic
      fibration $g:X\to\P^1$
      and $\mult_A(R_g)=1$. The list above of possible
      cycles in the affine case leads to this classification
      in the hyperbolic case.
      $$\G=T_{2,4,5}: \ A=(3;2,2,1,1;2,2,1;1).$$
      $$\G=T_{2,3,7}:\ A=(5;4,4,3,3,2,1;3,2;2).$$
      Everything else is impossible. In particular, $\G$ cannot
      be $T_{3,3,4}$.

      Finally, since $A^2=-2$, if $\G=T_{2,3,6}$
      then $A=(5;4,4,3,2,1;3,2;2).$

      \begin{proposition} \label{cond descr}\label{3.11}
        Suppose that $X$ is an Enriques surface.
        
        \part[i] $X$ is exceptional
        if and only if the support $\G$ of its
        conductrix $A$ is not negative definite.

        \part[ii] $X$ is exceptional if and only if $\G$ is
        one of the five diagrams just listed
        and $A$ is as listed.
        \begin{proof} It is enough to note that
          if $A$ is one of these cycles then
          $2A$ contains a Kodaira--N{\'e}ron
          fibre $F$, and $h^1(F,\sO_F)=1.$
        \end{proof}
      \end{proposition}

      We refer to this diagram $T_{p,q,r}$, or
      the triple $(p,q,r)$, as the \emph{type}
      of an exceptional Enriques surface.

      \begin{theorem}
        An exceptional Enriques
        surface $X$
        has a unique elliptic fibration if its
        type is affine. In this case the support of its conductrix is a half-fibre.
        $X$ has no elliptic fibration if its
        type is hyperbolic.
        \begin{proof} Immediate.
        \end{proof}
      \end{theorem}

      Recall that a genus $1$ fibration $f:X\to\P^1$
      is \emph{special} if it has a $2$-section that
      is a $(-2)$-curve and that,
      according to \cite{CD89}, Th. 3.4.1,
      if $f$ has no special $2$-section then there
      is another genus $1$ fibration $g:X\to \P^1$.

      \begin{lemma} If $X$ is exceptional of type $(2,4,4)$
        or $(2,3,6)$ then the unique elliptic
        fibration $f:X\to\P^1$ is not special.
        \begin{proof} Suppose that $C$ is a special $2$-section
          of $f$. Then $C.A>0$, so that, since $C$ is not
          contained in $A$, $C.A=1$. Moreover, $C$ meets
          $A$ in a point on a component $D$ of multiplicity
          $1$ in the Kodaira--N{\'e}ron half-fibre $F$
          supported on $A$. So $D$ is at the end of one of the
          long arms of the diagram. However, $A.D=0$
          in each case, and we have a contradiction
          to Lemma \ref{3.2}.
        \end{proof}
      \end{lemma}

      \begin{theorem}\label{classif}\label{3.14}
        If $X$ is exceptional then it is
        of type $(3,3,3),\ (2,4,5)$ or $(2,3,7)$.
        \begin{proof} It remains to exclude the types
          $(2,4,4)$ and $(2,3,6)$.

          If $X$ is of type $(2,4,4)$
        there is a unique elliptic fibration $f:X\to\P^1$
        and a quasi--elliptic $1$ fibration $g:X\to\P^1$.
        The curve $R_g$ has multiplicity $1$ in $A$
        and $A-R_g$ is $g$-vertical, so there is a fibre
        $X_t$ of $g$ that contains a cycle $D$
        of type either $A_7$ or $E_7$ and $A=A'+R_g$
        where $A'$ is supported on $D$.
        \begin{enumerate}
          \item Suppose $D_{red}$ is of type $A_7$.
          Then $A'=(1,2,2,3,2,2,1)$ and there is a component
          $C$ of $X_t$ that meets $A'$. Necessarily $C$ meets
          $A'$ in an end component $E$, so $A.C\ge 1$
          while $A.E=0$, contradiction.

        \item Suppose $D_{red}$ is of type $E_7$.
          Then $A'=(3;2,2,1;2,2;1)$ and we again get
          a contradiction by considering a component of
          $X_t$ that meets $A'$.
        \end{enumerate}

        If $X$ is of type $(2,3,6)$ then
        $A=A'+R_g$ and $A'$ has no component of
        multiplicity $1$, while there must be a $(-2)$-curve
        $C$ in the fibre containing $A'_{red}$ that meets $A'$,
        contradiction.
    \end{proof}
  \end{theorem}

  \begin{theorem}\label{3.15} Suppose that $X$ is exceptional of type
    $(p,q,r)$ and that its conductrix is $A$.
    
    \part[i] If $(p,q,r)=(3,3,3)$ then $A$ is supported on
    a half-fibre, of type $\tE_6$, of the unique elliptic fibration on $X$.

    \part[ii] If $(p,q,r)=(2,3,7)$ then $X$ has a unique genus $1$
    fibration $g:X\to \P^1$, $g$ is quasi--elliptic and
    $A-R_g$ is supported on a half-fibre of $g$ of type $\tE_8$.

    \part[iii] If $(p,q,r)=(2,4,5)$ then $X$ has no elliptic
    fibration and there is a unique quasi--elliptic fibration
    $g:X\to\P^1$ such that $A-R_g$ is supported on the whole
    of a fibre. This fibre is a half-fibre and is of type $\tE_7$.
    \begin{proof} This
      follows at once from the description
      of $A$ and the facts that $\mult_A(R_g)=1$ and
      $A-R_g$ is $g$-vertical.
    \end{proof}
  \end{theorem}

  \begin{proposition}\label{Aiii}
    If $X$ has a quasi--elliptic fibration
          $g:X\to\P^1$
          with a fibre $X_s$ (simple or double) of type $\tE_{7,8}$ then $X$ is
          exceptional. If the fibre is simple then $X$
          is of type $(3,3,3)$ or $(2,4,5)$
          accordingly.
          \begin{proof} $R_g$ has multiplicity $1$ in $A$ and
            $R_g.X_s=1$ or $2$. Moreover, $(R_g\cap X_s)_{red}$ is a single
            point, so $R_g$ meets $X_s$ in a curve of multiplicity
            $1$ or $2$ in $X_s$.

            Assume that $X$ is not exceptional, so that $\Supp A$ is
            negative definite. Then there are components $C_1,C_2$ of
            $X_s$ that do not lie in $A$ while $C_2$ meets $A$
            and $C_1.C_2=1$. Then $A.C_2\ge 1$ and $A.C_1\ge 0$,
            which contradicts Lemma \ref{transverse}. So
            $X$ is exceptional.

            If $X_s$ is simple then $R_g$ meets $X_s$ transversely
            in a component $C$ of multiplicity $2$ in $X_s$
            and the result follows from inspection of
            the possibilities provided by Proposition \ref{3.11}
            and Theorem \ref{3.14}.
          \end{proof}
        \end{proposition}
\end{section}
\begin{section}{From configuration to conductrix}
  In this section we describe exceptional Enriques surfaces $X$
  in terms of the configurations of $(-2)$-curves on them.

  \begin{lemma}\label{D4} If $f:X\to\P^1$ is an
    elliptic fibration and $f^{-1}(s)$ is a fibre
    that contains exactly one irreducible
    component $E$ of the conductrix $A$,
    then $f^{-1}(s)$ is of type $\tD_4$,
    $E$ is the branch vertex of the configuration
    and $\mult_A(E)=1$.
    \begin{proof}
      This is a consequence of Lemma \ref{3.2}.
    \end{proof}
  \end{lemma}
  
  \begin{theorem}\label{descript}
    An Enriques surface $X$ is exceptional if and only
    if it has a special genus $1$ fibration $g:X\to\P^1$ with a double
    fibre of type $\tE_{6,7,8}$. The type of $X$ is $(3,3,3)$,
    $(2,4,5)$ or $(2,3,7)$, accordingly.
    \begin{proof} Assume that $X$ is exceptional. We consider the
      various types separately.

      $T_{3,3,3}$: there is an elliptic fibration $f:X\to\P^1$
      with a double fibre $X_{f,s}$ of type $\tE_6$ and $A=(2;1,1;1,1;1,1)$.
      If $f$ is special we are done, and if not then there is,
      by Th. 3.4.1 of \cite{CD89}, another  genus $1$ fibration
      $g:X\to\P^1$. This is quasi--elliptic. Write $A=A'+R_g$, so
      that $A'$ is $g$-vertical. So there is a singular fibre
      $g^{-1}(t)=d_tX_{g,t}$ such that
      $\Supp A'\subsetneq \Supp X_{g,t}$.

      If $R_g$ is not an end curve of $A$ then
      $A'$ is disconnected; say $A'=A'_1\cup A'_2$,
      where $A'_2$ is an end curve of $A$. The
      fibre of $g$ that contains $A'_2$ is of type $\tD_4$,
      by Lemma \ref{D4},
      so has $5$ components. The fibre containing $A'_1$
      has at least $6$ components. But $\rho(X)=10$,
      contradiction.

      So $R_g$ is an end curve of $A$, so that $A'$ is of type $E_6$
      and $X_{g,t}=\tE_{r\ge 6}$. If $r=6$ let $E$ denote the curve such that
      $X_{g,t}=A'\cup E$, set-theoretically, and $C$ the curve in $A'$
      that meets $E$. Then $A.E\ge 1$
      and $A.C=1$, contradiction.

      So $r=7,8$ and $R_g$ is special.

      $T_{2,4,5}$ and $T_{2,3,7}$: the result is clear in these cases.

      Conversely, suppose that there is a special genus $1$
      fibration $g:X\to\P^1$ with a special $2$-section
      $D$ and a double fibre $g^{-1}(s)=2X_s$
      of type $\tE_r$. Suppose that $\dim(A\cap X_s)\le 0$.
      Then $c_1(\sF).E_i\ge 0$ for each component $E_i$
      of $X_s$, so that $c_1(\sF).E_i=0$ for all $E_i$
      and each $E_i$ is generically tangent to $\sF$.
      Let $E_1$ be the curve corresponding to the branch
      vertex of $\tE_r$, meeting $E_2,E_3,E_4$. Then
      $E_1\cap E_j$ is a singular point of $\sF$ for each $j=2,3,4$,
      but $c_1(\sF).E_1=0$, contradiction.
      So $A':= A\cap X_s$ is a non-zero divisor.
      \begin{enumerate}
        
      \item $r=6$. Then $X_s+D$ is a $T_{3,3,4}$
      configuration. If $g$ is elliptic then
      from Lemmas \ref{lem 1} to \ref{lem 3}
      it follows that $A=(2;1,1;1,1;1,1)$, and $X$ is exceptional.
      If $g$ is quasi--elliptic then $R_g\le A$ and $\G=T_{3,3,4}$
      which we know to be impossible.
      
    \item $r=7$. If $g$ is elliptic then $\mult_A(D)=0$
      and Lemmas \ref{lem 1} to \ref{lem 3} give a contradiction.
      If $g$ is quasi--elliptic then $\mult_A(D)=1$
      from Lemmas \ref{lem 1} to \ref{lem 3} it follows
      that $A=(3;2,2,1,1;2,2,1;1)$
      and $X$ is exceptional.
      
    \item $r=8$. As when $r=7$ we see that $g$ is quasi--elliptic,
      $$A=(5;4,4,3,3,2,1;3,2;2)$$
      and $X$ is exceptional.
    \end{enumerate}
  \end{proof}
\end{theorem}

\begin{corollary}\label{type}
  An Enriques surface $X$ is exceptional if and only
  if it contains a configuration of $(-2)$-curves of type $T_{p,q,r}$
  where $(p,q,r)=(3,3,4),(2,4,5)$ or $(2,3,7)$. The type of the
  surface is the type of the
  configuration except that a configuration
  $T_{3,3,4}$ gives a surface of type
  $T_{3,3,3}$.
  \begin{proof} If $X$ is exceptional then examination of the
    genus $1$ fibration provided by Theorem \ref{descript}
    gives the result.
    Conversely, a $T_{3,3,4}$ configuration
    yields a special genus $1$ fibration
    with a double fibre of type $\tE_6$, while
    $T_{2,4,5}$, resp., $T_{2,3,7}$, gives a special genus $1$
    fibration
    with a double fibre of type $\tE_7$, resp., $\tE_8$.
  \end{proof}
\end{corollary}

\begin{theorem}\label{4.4} Suppose that $X$ is an exceptional
  Enriques surface of type $T$ and conductrix $A$.

  \part[i] If $T=(2,3,7)$ then the only
  $(-2)$-curves on $X$ are the ones in $A$.

  \part[ii] If $T=(2,4,5)$ then
  there are exactly two $(-2)$-curves that are not in $A$.
  They form a fibre
  of type $Ta$ in the natural quasi--elliptic fibration
  $g:X\to\P^1$ given by Theorem \ref{3.15}. If this fibre
  has multiplicity $d$ then $X$ possesses just $3-d$ further
  genus $1$ fibrations. Each is quasi--elliptic and has a simple fibre
  of type $\tE_8$.
  \begin{proof} \DHrefpart{i}
    $E$ must meet $A$ since $\rho(X)=10$.
    If $E$ is not in $A$ then $A.E=1$, so $E$ meets $A$ in
    its unique component $R$ of multiplicity $1$. Since $R=R_f$
    where $f:X\to\P^1$ is the unique genus $1$ fibration
    on $X$ and $E$ is disjoint from $A-R$, $E$ is $f$-vertical.
    But this contradicts $\rho(X)=10$.

    \DHrefpart{ii} There is a further reducible fibre
    $f^{-1}(s)=dX_s$ of $f$; it has
    two components $E_1,E_2$ and $d=1$ or $2$.
    Since $X_s$ meets $R$ it meets $A$
    and since $\rho(X)=10$ it must be of type $Ta$.
    If $d=2$ then $R$ meets just one component
    and if $d=1$ then $R$ meets both components
    in their point of
    intersection. Examination of the diagram
    defined by $A\cup(E_1,E_2)$ concludes the proof.
  \end{proof}
\end{theorem}
Surfaces of type $(3,3,3)$ require more work and
are dealt with in the next section.
\end{section}
\begin{section}{$(3,3,3)$ exceptional surfaces}
  Suppose that $X$ is an exceptional Enriques surface
  of type $(3,3,3)$. We know that $X$ has a unique elliptic
  fibration $f:X\to\P^1$, that $f$ has a double fibre of type
  $\tE_6$ that supports the conductrix $A$ and that
  $A=(2;1,1;1,1;1,1)$.
  
  \begin{lemma} There exists at least one quasi--elliptic
    fibration on $X$.
    \begin{proof} If $f$ were the only one genus $1$ fibration
      on $X$ then, by Th. 3.4.1 of \cite{CD89}, it would
      have a fibre of type $\tE_8$. But $\rho(X)=10$
      so this is impossible.
    \end{proof}
  \end{lemma}

  \begin{lemma}\label{E7}
    \part[i] For any quasi--elliptic fibration
    $g:X\to\P^1$ the curve $R_g$ is an end curve of $A$
    and the fibre $X_s$ of $g$ that contains $A-R_g$ is
    a simple fibre of type $\tE_7$.

    \part[ii] Every $\tE_7$-configuration $\Psi$ of
    $(-2)$-curves on $X$ arises in this way.
    \begin{proof} \DHrefpart{i} See the proof of Theorem \ref{descript}
      for the fact that $R_g$ is an end curve of $A$.

      Say $g^{-1}(s)=d_sX_s$. Since $X_s$ contains the $E_6$
      configuration $A-R_g$ it is of type $\tE_r$. Lemma \ref{3.2}
      excludes the cases $r=6,8$ and shows that $R_g$ meets the short
      arm of a $\tE_7$ fibre transversely. So $d_s=1$.

      \DHrefpart{ii} There is a genus one fibration $g:X\to\P^1$
      such that $e\Psi$ is a fibre and $e=1$ or $2$. $T_0$ is not
      $h$-vertical, since $\rho(X)=10$, so that $g\ne f$ and therefore
      $g$ is quasi--elliptic. $R_g$ is an end curve of $A$ and
      $A-R_g$ is $g$ -vertical; then $A-R_g\subset\Psi$
      since again $\rho(X)=10$ and \DHrefpart{ii} is proved.
    \end{proof}
  \end{lemma}

  Write $\Supp A=T_0$.

  \begin{lemma}\label{unique}
    $T_0$ is the unique $T_{3,3,3}$-configuration
    on $X$.
    \begin{proof}
      Suppose $S$ is another such configuration.
      Then there is a quasi--elliptic fibration $g:X\to\P^1$
      with a fibre supported on $S$. The curve $R_g$ is an
      end component of $A$ and $A-R_g$ is $g$-vertical.
      Since $S$ and $A-R_g$ are connected and $\rho(X)=10$,
      it follows that $T_0-R_g\subset S$. Then there is an end
      curve $C$ of $S$ that is not contained in
      $A$ and which meets $A$ in a curve
      $D$ of multiplicity at least $2$. Therefore $C.A\ge 2$,
      which is impossible.
    \end{proof}
  \end{lemma}

  Let $\sS$ denote the set of $T_{4,4,4}$ configurations on $X$.

  \begin{lemma} $T_0$ extends to some $S\in\sS$ 
    and every element of $\sS$ 
    contains $T_0$.
    \begin{proof} Take the $\tE_7$ configuration that is the
      fibre provided by the fibre $X_s$ of
      Lemma \ref{E7}. There is another reducible
      fibre of $g$, since $\rho(X)=10$; it is of type $Ta$.
      A suitable component of this fibre extends $X_s+R$
      to a $T_{4,4,4}$ configuration.
      The rest follows from Lemma \ref{unique}.
    \end{proof}
  \end{lemma}
  
  Let $\rho:X\to Y_1=X/\sF$ be the quotient, so that
  $Y_1=\tZ^{(1)}$. We know that $Y_1$is a rational surface
  and that $\Sing Y_1=4\times A_1$ or
  $1\times D_4$.

  \begin{lemma}\label{smooth} If $C$ is a curve in $X$ such that
    $\sF\vert_C$ maps isomorphically to either $T_C$ or $N_{C/X}$
    then $Y_1$ is smooth along $D=\rho(C)$. If $\sF\vert_C=T_C$
    then $\rho_*C=2D,\rho^*D=C$ and $D^2=C^2/2$.
     If $\sF\vert_C=N_{C/X}$
     then $\rho_*C=D,\rho^*D=2C$ and $D^2=2C^2$.
     \begin{proof} This is standard.
     \end{proof}
   \end{lemma}

   \begin{proposition}
     \part[i] $T_0$ maps to a normally crossing
     configuration $U_1$ of $\P^1$'s on $Y_1$,
     disjoint from $\Sing Y_1$ and described by
     $U_1=(-4;-1,-4;-1,-4;-1,-4)$.

     \part[ii] There is a birational
     contraction $\pi:Y_1\to Y$ of the central $D_4$
     configuration in $U_1$ to a smooth point. The
     image of $U_1$ is a Kodaira--N{\'e}ron cycle
     $U_0=C_1+C_2+C_3$ of type $Tr$. There is 
     a commutative diagram
     $$
     \xymatrix
     {
       {X}\ar[r]\ar[rd]_f&{Y_1}\ar[r]^{\pi}\ar[d]^{h_1}&{Y}\ar[dl]^{h}\\
       &{\P^1}&
     }
     $$

     \part[iii] $h:Y\to\P^1$ is a relatively minimal
     elliptic fibration and $U_0\in\vert -K_Y\vert$.

     \part[iv] The quasi--elliptic fibrations on $X$ correspond
     to the rulings on $Y$.

     \part[v] Give a ruling $q:Y\to\P^1$, two of the curves $C_i$
     are $q$-vertical and the third is a purely inseparable
     $2$-section of $q$.

     \part[vi]\label{type H} Every $S\in\sS$
     maps to a configuration $H=U_0+\sum_1^3D_j$
     of $\P^1$'s on $Y-\Sing(Y)$ such that
     $$C_i.D_j=-D_i.D_j=\delta_{ij}.$$

     \part[vii] Write $\Pic Y=L$. Then $L$ is of signature $(1,5)$ and
     $L^\perp/L\cong (\Z/2)^2$.

     \part[viii] $\MW(X)$ is the Mordell--Weil rank of $h:Y\to\P^1$.
     \begin{proof} Where this does not
       follow from what we already know it is easy.
     \end{proof}
   \end{proposition}
   
   A \emph{type $H$ configuration of curves} on $Y$
   is a configuration $U_0+\sum_1^3D_j$ of irreducible curves
   as in \ref{type H}. A \emph{type $H$ configuration
     of classes} on $Y$ is the same thing, except that each
   $D_j$ is only required to be a class in $\Pic(Y)$. So $\sS$
   is identified with the set of type $H$ configurations
   of curves on $Y$.
   
   \begin{lemma} Suppose that $E$ is a section of $h:Y\to\P^1$.
     If $E^2<0$ then $Y$ is smooth along 
     $E$ and $E^2=-1$. Conversely, if $Y$ is smooth along $E$
     then $E^2=-1$.
     \begin{proof} $E$ is disjoint from the exceptional locus
       of $Y_1\to Y$, so that we can write $E=\rho(C)$
       for some curve $C$ on $X$ with $C^2<0$.
       Therefore $C$ is a $(-2)$-curve. Since $C$ meets $A$,
       $\sF\vert_C$ maps isomorphically to $T_C$, and we can apply
       Lemma \ref{smooth}. The converse is well known.
     \end{proof}
   \end{lemma}

   \begin{lemma}\label{section} Suppose that $V$ is an RDP surface,
     that $q:V\to \P^1$ is an elliptic fibration with fibre $\Phi$
     and that $D$ is a Cartier divisor class on $V$ such that
     $D.\Phi=1$ and that $D^2=-\chi(\sO_V)$. Assume also that
     $D$ and $K_V$ are nef relative to $q$.

     Then $D$ is the class of a section of $q$
     that is disjoint from $\Sing V$.
     \begin{proof} Since $D$ is Cartier, $D.\Phi=1$ and
       $D$ is $q$-nef it follows that there is a birational
       contraction
       $\tau:V\to V'$ that contracts exactly the $q$-vertical
       curves
       $\psi$ with $D.\psi=0$, $V'$ also has RDPs
       and $D=\tau^*D'$ for some Cartier
       divisor class $D'$ on $V'$. Then $q$ factors as
       $q=q'\circ\tau$ with $q':V'\to \P^1$
       and all fibres
       of $q':V'\to\P^1$ are reduced and irreducible.
       Let $\Phi'$ be a fibre of $q'$. Then
       there is a divisor class $a$ on $C$ such that
       $D'+q'^*a$ is effective; if $a$ is taken to have minimal
       degree then $D'+q'^*a\sim D_1$ where
       $D_1$ is a section. Then
       $$-\chi(\sO_{V'})=D_1^2=(D')^2+2\deg a=-\chi(\sO_V)+2\deg a,$$
       so that $a=0$ and $D'$ is the class of a section.
       So $D\sim\tau^*D'_1$ and is therefore the class of a section.
       Since this section is a Cartier divisor on $V$ and is smooth it is
       disjoint from $\Sing V$.
     \end{proof}
   \end{lemma}

   A $(-2)$-curve $E$ on $X$ is \emph{extraneous} if it is disjoint from
   $T_0$. Equivalently, $E$ is extraneous if it is an
   irreducible component of a reducible fibre of $f$ besides $T_0$.
   $E$ is \emph{horizontal} if it is not $f$-vertical, or,
   equivalently,
   if it is not in $T_0$ and is not extraneous.
   A curve on $Y$ is extraneous if it is the image of an extraneous
   curve on $X$. A fibre of $f$ or $h$ is extraneous if it consists
   of extraneous curves. Extraneous curves (or fibres) exist
   if and only if $\MW(X)\le 1$.
   
   \begin{proposition} Assume
     that $\MW(X)=2$, that $D\in L$, $D^2=-1$,
     $D.C_i\ge 0$ for all $i$ and $D.U_0=1$. Then $D$ is the class of a section
     of $h:Y\to\P^1$.
     \begin{proof} $h$ has no extraneous fibres
       and the result follows from Lemma \ref{section}.
     \end{proof}
   \end{proposition}

   Put $M=\sum\Z.C_i\subset L$ and
   $\Delta=\{\g\in O(L)\vert \g(C_i)=C_i\ \forall\ i\}$.
   Fix an $H$-configuration $U_0\cup\{D_j\}$ on $Y$
   and put $\phi_i=D_j+C_j+C_k+D_k$ when
   $\{i,j,k\}=\{1,2,3\}$. Put
   $\a_i=D_i-\phi_i/2\in M^\perp\subset L^\vee$,
   so that $\a_i.\a_j=(1-3\delta_{ij})/2$.
   Let $s_i$ be the reflexion in $\a_i$; then
   $s_i\in\Delta$. Define
   $W=\langle s_1,s_2,s_3\rangle\subset\Delta$.

   \begin{lemma}\label{weyl}
     
     \part[i] $W$ acts on $M^\perp$ as the
     Weyl group $W(\tA_2)$.

     \part[ii] $\Delta=W$.
     \begin{proof} \DHrefpart{i} $s_is_j$ has order
       $3-\delta_{ij}$. So there is a surjection $\pi:W(\tA_2)\to W$.
       But the reflexion group action of $W(\tA_2)$ visibly
       factors through $\pi$.

       \DHrefpart{ii} This follows from the facts that
       $W(\tA_2)\times (\pm 1)$ is the full group
       $O(M^\perp)$ and that $\Delta$ acts effectively
       on $M^\perp$.
     \end{proof}
   \end{lemma}
   
   \begin{theorem}\label{MW2} Suppose that $\MW(X)=2$.

     \part[i] $\sS$ is a torsor under $W(\tA_2)$.

     \part[ii] Each horizontal $(-2)$-curve on $X$
     lies in exactly six elements of $\sS$.

     \part[iii] Each $S\in \sS$ defines three quasi--elliptic fibrations
     on $X$ and each quasi--elliptic fibration on $X$
     arises from two such diagrams.
     \begin{proof} \DHrefpart{i} There are no extraneous curves
       and therefore any $(-2)$-curve on $X$ 
       corresponds to a $(-1)$-curve on $Y$
       that meets $U_0$. We know that every type $H$
       configuration of classes on $Y$ is
       given by a unique type $H$ configuration of curves on $Y$,
       so that $\sS$ is identified with the set of type
       $H$ configurations of classes on $Y$. Therefore
       $\sS$ is a torsor under $\Delta$.
       \DHrefpart{i} now follows from this and Lemma \ref{weyl}.

       For \DHrefpart{ii} and \DHrefpart{iii} regard $W$ as a symmetry
       group of a tessellation of the Euclidean plane
       $\Pi=(M^\perp/M\cap M^\perp)\otimes\R$
       by equilateral triangles each of which is a fundamental
       domain for the action of $W$ on $\Pi$. \DHrefpart{ii} follows from
       noting that each vertex of the tessellation
       lies in six triangles while \DHrefpart{iii} follows from the
       facts that each triangle has three edges and each edge lies in
       two triangles.
     \end{proof}
   \end{theorem}

   \begin{lemma} Suppose that $\MW(X)=1$.

     \part[i] $h$ has one extraneous fibre $\Phi_0$, say,
     with components $\psi_1,\psi_2$.

     \part[ii] Suppose that $D\in L$ with $D^2=-1$, $D.U_0=1$
     and that $D.\psi_i\ge 0$.
     Then $D$ is the class of a section of $h:Y\to\P^1$.
     \begin{proof} \DHrefpart{i} is clear. For \DHrefpart{ii},
       apply Lemma \ref{section}.
     \end{proof}
   \end{lemma}

   \begin{theorem}\label{MW1} Suppose that $\MW(X)=1$.

     \part[i] $\sS$ is a torsor
     under the infinite cyclic subgroup of $W$
     generated by some glide-reflexion.

     \part[ii] Each horizontal $(-2)$-curve on $X$
     lies in exactly three elements of $\sS$.

     \part[iii]  Each $S\in \sS$ defines three quasi--elliptic fibrations
     on $X$ and each quasi--elliptic fibration on $X$
     arises from two such diagrams.
     \begin{proof} We have to classify type $H$
       diagrams of classes on $Y$ such that
       each $D_i$ is $h$-nef. That is,
       $D_i.\psi_j\ge 0$ for all $i,j$. Since each type $H$
       diagram spans $L$, we require, after renumbering
       if necessary, that $D_1.\psi_1=D_2.\psi_2=D_3.\psi_2=1$
       and $D_i.\psi_j=0$ otherwise. Again, regard $W$ as a group
       of symmetries of 
       $\Pi$ that preserves
       a tessellation into equilateral triangles; one sees that the
       theorem follows and that $\g$ takes one triangle
       to another with a common edge.
       So $\langle\g\rangle$ acts on a strip $\Sigma$
       whose width is one triangle.
       A single triangle forms a fundamental domain
       for the action of $\langle\g\rangle$ on $\Sigma$.
       Since each vertex of the tessellation of $\Sigma$
       lies in three triangles in $\Sigma$ \DHrefpart{i} is proved.

       \DHrefpart{ii} and \DHrefpart{iii} are proved as in
       Theorem \ref{MW2}.
     \end{proof}
   \end{theorem}

   \begin{lemma}\label{one fibre} If $\MW(X)=0$ then $h:Y\to\P^1$
     has just one extraneous fibre $\Phi_0$.
     It is simple and has three components.
     \begin{proof} It is enough to prove the
       analogous result for $f:X\to\P^1$. Recall
       that $T_0$ is a double fibre of $f$.
       Suppose that $v_1,v_2,v_3$ are the end curves
       in a $T_{4,4,4}$ diagram. Then $v_i$ is a special bisection of $f$.
       Suppose that $\Phi_0, \Phi_1$ are reducible fibres
       each of which has
       just two components, say $\psi_1,\psi_2$ and $\psi_3,\psi_4$,
       respectively.
       Then $A.\psi_j=0$ and $A.v_i=1$, so that $v_i$ cannot
       meet $\psi_j$ transversely. So $\Phi_k$ is a simple fibre,
       $v_i.\psi_j=0$ or $2$.

       After renumbering if necessary, $\psi_1.v_i=2\delta_{1i}$
       and $\psi_3.v_i=2\delta_{3i}$. Then
       $$(v_1+\psi_1)^2=0=(v_3+\psi_3)^2=(v_1+\psi_1).(v_3+\psi_3),$$
       so that $v_1+\psi_1$
       and $v_3+\psi_3$ are proportional. However, if $u_3$
       is the end curve in $T_0$ that meets $v_3$ then
       $u_3.(v_3+\psi_3)=1$ while $u_3.(v_1+\psi_1)=0$.
       This contradiction proves the lemma.
     \end{proof}
   \end{lemma}
   
   \begin{theorem}\label{MW0} If $\MW(X)=0$ then $\sS$ has one element
     and there are only three horizontal $(-2)$-curves
     on $X$.
     \begin{proof} As before, we must classify
       type $H$ diagrams of classes on $Y$ where
       each $D_i$ is $h$-nef. By Lemma \ref{one fibre}
       $h$ has just one reducible fibre besides $U_0$,
       with three components $\psi_i$. We can take
       $D_i.\psi_j=\delta_{ij}$ and now the result is immediate.
     \end{proof}
   \end{theorem}
\end{section}

\end{document}